\documentclass[a4paper]{article}

\usepackage{amsfonts, latexsym, amssymb, amsthm, amsmath, mathrsfs, esint}

\newcommand{\R}{\mathbb{R}}
\newcommand{\N}{\mathbb{N}}

\newcommand{\Ha}{\mathcal{H}}

\newcommand{\loc}{\mathrm{loc}}

\let\div\relax\DeclareMathOperator{\div}{div}
\DeclareMathOperator{\curl}{curl}

\DeclareMathOperator{\dist}{dist}

\newcommand{\blank}{{\mkern 2mu\cdot\mkern 2mu}}

\newcommand{\dd}[2]{\frac{\partial #1}{\partial #2}}
\newcommand{\set}[2]{\left\{ #1 \colon #2 \right\}}

\newcommand{\restr}{\mathchoice
{\kern2pt\mbox{\vrule width 0.08ex height1.5ex depth0ex\kern-0.08ex\vrule width 1.5ex height.08ex depth0ex}\kern2pt}
{\kern2pt\mbox{\vrule width 0.08ex height1.5ex depth0ex\kern-0.08ex\vrule width 1.5ex height.08ex depth0ex}\kern2pt}
{\kern1.5pt\mbox{\vrule width 0.06ex height1.1ex depth0ex\kern-0.06ex\vrule width 1.1ex height.06ex depth0ex}\kern1.5pt}
{\kern1pt\mbox{\vrule width 0.04ex height0.75ex depth0ex\kern-0.04ex\vrule width 0.75ex height.04ex depth0ex}\kern1pt}
}

\newtheorem{theorem}{Theorem}
\newtheorem{lemma}[theorem]{Lemma}
\newtheorem{proposition}[theorem]{Proposition}
\newtheorem{definition}[theorem]{Definition}

\theoremstyle{remark}
{}
\newtheorem*{remark}{Remark}

\begin{document}

\title{The streamlines of $\infty$-harmonic functions obey the inverse mean curvature flow}

\author{Roger Moser\footnote{Department of Mathematical Sciences,
University of Bath, Bath BA2 7AY, UK.
E-mail: r.moser@bath.ac.uk}}

\maketitle

\begin{abstract}
Given an $\infty$-harmonic function $u_\infty$ on a domain
$\Omega \subseteq \R^2$, consider the function $w = -\log |\nabla u_\infty|$.
If $u_\infty \in C^2(\Omega)$ with $\nabla u_\infty \neq 0$ and
$\nabla |\nabla u_\infty| \neq 0$, then it is easy to check that
\begin{itemize}
\item the streamlines of $u_\infty$ are the level sets of $w$ and
\item $w$ solves the level set formulation of the inverse mean curvature flow.
\end{itemize}

For less regular solutions, neither statement is true in general,
but even so, $w$ is still a weak solution of the inverse
mean curvature flow under far weaker assumptions. This is proved
through an approximation of $u_\infty$ by $p$-harmonic functions,
the use of conjugate $p'$-harmonic functions, and the
known connection of the latter with the inverse mean curvature flow.
A statement about the regularity of $|\nabla u_\infty|$ arises as a by-product.
\end{abstract}

\section{Introduction}

Let $\Omega \subseteq \R^2$ be an open set. A function $u_\infty \in C^0(\Omega)$
is called $\infty$-harmonic if it is a viscosity solution of the Aronsson
equation
\begin{equation} \label{eqn:infty-harmonic}
\left(\dd{u_\infty}{x_1}\right)^2 \frac{\partial^2 u_\infty}{\partial x_1^2} + 2\dd{u_\infty}{x_1} \dd{u_\infty}{x_2} \frac{\partial^2 u_\infty}{\partial x_1 \partial x_2} + \left(\dd{u_\infty}{x_2}\right)^2 \frac{\partial^2 u_\infty}{\partial x_2^2} = 0.
\end{equation}
This equation was introduced by Aronsson \cite{Aronsson:67, Aronsson:68}, motivated by optimal
Lipschitz extensions of the boundary data, and has been studied extensively
since then. Highlights of the theory include existence \cite{Bhattacharya-DiBenedetto-Manfredi:89} and uniqueness \cite{Jensen:93} of solutions for
boundary value problems associated to \eqref{eqn:infty-harmonic},
regularity results \cite{Savin:05, Evans-Savin:08}, and connections to
stochastic tug-of-war games \cite{Peres-Schramm-Sheffield-Wilson:09}.

For $u_\infty \in C^2(\Omega)$, equation \eqref{eqn:infty-harmonic} may alternatively be represented as
\[
\nabla u_\infty \cdot \nabla |\nabla u_\infty|^2 = 0.
\]
It is then obvious that the function $|\nabla u_\infty|$ is constant
along the streamlines of $u_\infty$, i.e., along the curves in
$\Omega$ arising through the solutions of the ordinary differential
equation $\dot{\gamma}(t) = \nabla u_\infty(\gamma(t))$.
This is one of the reasons why the streamlines of an $\infty$-harmonic function are of
particular interest and have received some attention in the literature
\cite{Aronsson:68, Lindgren-Lindqvist:19, Lindgren-Lindqvist:21}. In general,
however, viscosity solutions of \eqref{eqn:infty-harmonic} are not
$C^2$-regular. It was shown by Evans and Savin \cite{Evans-Savin:08} that they
are of class $C^{1, \alpha}$ for some $\alpha > 0$. As the example
$u_\infty(x) = |x_1|^{4/3} - |x_2|^{4/3}$ of Aronsson \cite{Aronsson:84}
shows, no exponent better than $\alpha = 1/3$ can be expected. Nevertheless,
at least in an annular domain with boundary values $0$ and $1$, respectively,
on the two boundary components, it was shown by Lindgren and Lindqvist
\cite{Lindgren-Lindqvist:19, Lindgren-Lindqvist:21} that
$|\nabla u_\infty|$ is constant along streamlines that are generic in some sense.
There is also a weaker statement that is true in general (see, e.g., the description by Crandall \cite[Section 6]{Crandall:08}).

If $u_\infty \in C^2(\Omega)$ is a solution of \eqref{eqn:infty-harmonic}
with $\nabla u_\infty \neq 0$ and $\nabla |\nabla u_\infty| \neq 0$,
then we also conclude that
\[
\frac{\nabla |\nabla u_\infty|}{|\nabla |\nabla u_\infty||} = \pm \frac{\nabla^\perp u_\infty}{|\nabla u_\infty|},
\]
where we write $\nabla^\perp = (-\dd{}{x_2}, \dd{}{x_1})$.
Hence the function $w = -\log|\nabla u_\infty|$ satisfies
\[
\div\left(\frac{\nabla w}{|\nabla w|}\right) = \mp \div\left(\frac{\nabla^\perp u_\infty}{|\nabla u_\infty|}\right) = \pm \frac{\nabla^\perp u_\infty \cdot \nabla |\nabla u_\infty|}{|\nabla u_\infty|^2} = |\nabla w|.
\]
The equation
\begin{equation} \label{eqn:imcf}
\div\left(\frac{\nabla w}{|\nabla w|}\right) = |\nabla w|
\end{equation}
has a geometric interpretation: it is the level set formulation of the inverse mean curvature flow.

The inverse mean curvature flow is an evolution equation for
hypersurfaces. It is often studied on a Riemannian manifold, but we explain
it here for an open set $\Omega \subseteq \R^n$. Consider an oriented
$(n - 1)$-dimensional manifold $N$ and a
smooth map $\phi \colon [0, T) \times N \to \Omega$ such that
$N_t = \phi(\{t\} \times N)$ is an immersed hypersurface 
for every $t \in [0, T)$. Suppose that
$\nu \colon [0, T) \times N \to S^{n - 1}$ is a smooth map
such that $\nu(t, \blank)$ is a normal vector field on
$N_t \subseteq \Omega$ for every $t$,
and let $H \colon [0, T) \times N \to \R$ be the function such that
$H(t, \blank)$ is the corresponding (scalar) mean curvature of $N_t$.
We say that $\phi$ is a classical solution of the inverse mean curvature flow if
\begin{equation} \label{eqn:imcf-classical}
\dd{\phi}{t} = \frac{\nu}{H}
\end{equation}
in $(0, t) \times N$. This is a parabolic equation, so we may hope
to solve it for a prescribed initial hypersurface $N_0$ under suitable boundary
conditions. But this is
not always possible, either because $H$ has zeroes at $t = 0$ or because singularities develop in finite time. For this reason,
a weak notion of solutions was proposed by Huisken and Ilmanen
\cite{Huisken-Ilmanen:01}, based on a level set formulation.
The underlying idea is to look for a function $w \colon \Omega \to \R$
such that $N_t = w^{-1}(\{t\})$. As long as $w$ is sufficiently
smooth and $\nabla w \neq 0$, equation \eqref{eqn:imcf-classical}
is equivalent to \eqref{eqn:imcf}.

Equation \eqref{eqn:imcf}, however, allows a weak interpretation as well.
Huisken and Ilmanen use a variational principle for this purpose.
In this paper, we use a different formulation, which is
more convenient for our main results and perhaps more intuitive, too.
We will see in
Section \ref{sct:imcf}, however, that the following condition implies
that $w$ is a weak solution in the sense of Huisken and Ilmanen,
as long as we impose enough regularity such
that the latter makes sense.

\begin{definition} \label{def:imcf}
A function $w \in W^{1, 1}_\loc(\Omega)$ is called a \emph{weak solution
of \eqref{eqn:imcf}} if there exists a measurable vector field
$F \colon \Omega \to \R^n$ such that $|F| \le 1$ and
$F \cdot \nabla w = |\nabla w|$ almost everywhere in $\Omega$
and $\div F = |\nabla w|$ weakly in $\Omega$.
\end{definition}

We now restrict our attention to $n = 2$ again. For
$w = -\log |\nabla u_\infty|$, equation \eqref{eqn:imcf} means that
the level sets of $|\nabla u_\infty|$ move by the inverse mean
curvature flow. Furthermore, under the above regularity assumptions, the level
sets of $|\nabla u_\infty|$ are the streamlines of $u_\infty$.

We study the question to what extent these observations persist when
we remove the regularity assumptions. It is not true in general that
viscosity solutions of \eqref{eqn:infty-harmonic} give rise even to
weak solutions of \eqref{eqn:imcf}. The function
$u_\infty(x) = |x_1|^{4/3} - |x_2|^{4/3}$ provides a counterexample
here, too. In this case, as $\nabla|\nabla u_\infty| \neq 0$ almost
everywhere, there is only one possible choice for the  vector field $F$ from Definition \ref{def:imcf}. We can then check that \eqref{eqn:imcf}
does not hold on the coordinate axes. By contrast, the function
$u_\infty(x) = \xi \cdot x$, for a constant $\xi \in \R^2 \setminus \{0\}$,
may appear an unlikely candidate for the
inverse mean curvature flow. Its streamlines are straight lines
and the curvature vanishes identically. But in the formulation of Definition
\ref{def:imcf} (and also in the formulation of Huisken and Ilmanen
\cite{Huisken-Ilmanen:01}), the inverse mean curvature flow can
deal with this situation. It is easily seen that any constant function
is a weak solution of \eqref{eqn:imcf}. If $\nabla w = 0$,
the geometric interpretation is
that rather than moving continuously, the hypersurfaces jump
instantaneously. This is indeed the typical way in which the weak
inverse mean curvature flow resolves singularities.

We can prove that $w = -\log |\nabla u_\infty|$ does solve \eqref{eqn:imcf}
weakly under an additional assumption, which is based on the idea
that the function $|\nabla u_\infty|$ is monotone along the level sets of
$u_\infty$. This induces a sense of direction on the level sets
and, by comparison with $\nabla u_\infty$, an orientation of $\R^2$.
Writing $B_r(x)$ for the open disk in $\R^2$ with centre $x$ and radius $r$,
we can formalise this notion as follows.

\begin{definition}
Let $u_\infty \in C^1(\Omega)$ with $\nabla u_\infty \neq 0$ everywhere.
For $G \subseteq \Omega$, an \emph{orientation} of $u_\infty$ in $G$
is a continuous function $\omega \colon G \to \{-1, 1\}$
such that for any $x \in G$ there exists $r > 0$ with the
following property: for all $y, z \in B_r(x) \cap G$, if
$u_\infty(y) = u_\infty(z)$ and $(z - y) \cdot \nabla^\perp u_\infty(x) \ge 0$, then
\[
\omega(x) |\nabla u_\infty(z)| \ge \omega(x) |\nabla u_\infty(y)|. 
\]
\end{definition}

For example, the function $u_\infty(x) = |x_1|^{4/3} - |x_2|^{4/3}$
has no orientation in $\R^2$, but does have an orientation in
$\set{x \in \R^2}{x_1 x_2 \neq 0}$, which is
$\omega(x) = \frac{x_1 x_2}{|x_1||x_2|}$.
More generally, if $u_\infty \in C^2(\Omega)$ with $\nabla u_\infty \neq 0$
and $\nabla |\nabla u_\infty| \neq 0$, then $u_\infty$ has an orientation
in $\Omega$, which coincides with the orientation of $\R^2$
induced by the
pair of vectors $(\nabla u_\infty, \nabla |\nabla u_\infty|)$.
The function $\nabla u_\infty(x) = \xi \cdot x$ also has an orientation (in fact more than one).

\begin{theorem} \label{thm:imcf}
Let $u_\infty \in C^1(\Omega)$ be an $\infty$-harmonic function
with $\nabla u_\infty \neq 0$ in $\Omega$.
Suppose that $\omega$ is an orientation of $u_\infty$ in $\Omega$.
Then $|\nabla u_\infty|$
belongs to $W_\loc^{1, q}(\Omega)$ for all $q < \infty$ and satisfies
\begin{equation} \label{eqn:gradient}
|\nabla |\nabla u_\infty|| \frac{\nabla^\perp u_\infty}{|\nabla u_\infty|} = \omega \nabla |\nabla u_\infty|
\end{equation}
almost everywhere and
\begin{equation} \label{eqn:div}
\div \left(\frac{\nabla^\perp u_\infty}{|\nabla u_\infty|}\right) = - \omega \frac{|\nabla |\nabla u_\infty||}{|\nabla u_\infty|}
\end{equation}
weakly in $\Omega$. Hence the function
$w = - \log |\nabla u_\infty|$ is a weak solution of \eqref{eqn:imcf}.
\end{theorem}

Equation \eqref{eqn:gradient} may be regarded as another representation
of the Aronsson equation \eqref{eqn:infty-harmonic}.
Equation \eqref{eqn:div}, on the other hand, provides additional
information about the behaviour of the solutions.

It is already known
from work of Koch, Zhang, and Zhou \cite{Koch-Zhang-Zhou:19.1}
that $|\nabla u_\infty| \in W_\loc^{1, 2}(\Omega)$.
Theorem \ref{thm:imcf} improves this regularity to $W^{1, q}_\loc(\Omega)$
for any $q < \infty$,
but only under the additional conditions that $\nabla u_\infty \neq 0$
and there is an orientation. As discussed in the aforementioned
paper, this is false in general. A counterexample is given
by our usual suspect $u_\infty(x) = |x_1|^{4/3} - |x_2|^{4/3}$.
Nevertheless, the regularity statement from Theorem \ref{thm:imcf}
can be improved somewhat.

\begin{theorem} \label{thm:regularity}
Let $u_\infty \in C^1(\Omega)$ be an $\infty$-harmonic function
with $\nabla u_\infty \neq 0$ in $\Omega$. Let $G \subseteq \Omega$
be an open set and $\Gamma \subset \partial G \cap \Omega$. Suppose that
$\omega$ is an orientation of $u_\infty$ in $G \cup \Gamma$. Suppose further that
for every $x \in \Gamma$ there exist $r > 0$ and a Lipschitz
function $f \colon \R \to \R$ such that
\[
G \cap B_r(x) = \set{y \in B_r(x)}{\omega(x) y \cdot \nabla^\perp u_\infty(x) < f(y \cdot \nabla u_\infty(x))}.
\]
Let $U \subseteq G$ be an open set with $\overline{U} \subseteq G \cup \Gamma$.
Then $|\nabla u_\infty| \in W^{1, q}(U)$ for every $q < \infty$.
\end{theorem}

In less technical terms, we require that $|\nabla u_\infty|$ is non-decreasing
if we travel along a level set of $u_\infty$ inside $G$ towards $\Gamma$.
There is no such restriction outside of $\Gamma$.
In some cases, when $|\nabla u_\infty|$
has local maxima on $\Gamma$, it may be possible to apply the theorem
on the other side of $\Gamma$ as well with the opposite orientation.
Even then, however, it does not follow that $-\log |\nabla u_\infty|$ will
satisfy equation \eqref{eqn:imcf} on $\Gamma$.

The assumption that $f$ is Lipschitz continuous is stronger than
necessary; we use it for the sake of a simpler statement.
A weaker assumption is used in Proposition \ref{prp:local} below.

Clearly we need some prior information about the behaviour of
$u_\infty$ before we can apply a result such as this.
Such information is available, for example, for the $\infty$-harmonic
functions studied by Lindgren and Lindqvist
\cite{Lindgren-Lindqvist:19, Lindgren-Lindqvist:21}. Combining their
results with Theorem \ref{thm:regularity}, we see that under the
assumptions of the second paper \cite{Lindgren-Lindqvist:21},
we have local $W^{1, q}$-regularity of $|\nabla u_\infty|$ for all
$q < \infty$ away
from what Lindgren and Lindqvist call the attracting streamlines.

The main purpose of this paper, however, is not to provide regularity
results, but to explore the relationship between $\infty$-harmonic functions and the inverse
mean curvature flow. It seems that this has not been discussed
in the literature before even in the smooth case, although
some related calculations are present in the work of
Aronsson \cite{Aronsson:68} and Evans \cite{Evans:93}.
The proof of Theorem \ref{thm:imcf} shows that the connection is in fact deeper than
the simple calculations at the beginning of the introduction suggest.
The arguments are based on the following ideas, explained here for
$u_\infty \in C^1(\overline{\Omega})$ when $\Omega$ is a simply connected
domain with Lipschitz boundary.

According to the results of Jensen \cite{Jensen:93},
an $\infty$-harmonic function can be approximated by $p$-harmonic
functions. Let therefore $u_p$ denote the unique minimisers of the functionals
\[
E_p(u) = \left(\frac{1}{p} \int_\Omega |\nabla u|^p \, dx\right)^{1/p}
\]
in the spaces $u_\infty + W_0^{1, p}(\Omega)$. Then they
satisfy the $p$-Laplace equation
\[
\div(|\nabla u_p|^{p - 2} \nabla u_p) = 0.
\]
We use the notation $\Delta_p u = \div(|\nabla u|^{p - 2} \nabla u)$
for the $p$-Laplace operator; then we can write this equation in the form
$\Delta_p u_p = 0$.

We eventually consider the limit $p \to \infty$, but for the
moment we fix $p < \infty$. Let $p'$ be its conjugate exponent
with $\frac{1}{p} + \frac{1}{p'} = 1$. Then we
note that $\curl(|\nabla u_p|^{p - 2} \nabla^\perp u_p) = \Delta_p u_p = 0$.
Hence
there exists $v_p \in W^{1, p'}(\Omega)$ satisfying
\[
\nabla v_p = \omega |\nabla u_p|^{p - 2} \nabla^\perp u_p.
\]
Then we compute $|\nabla v_p|^{p' - 2} \nabla v_p = \omega \nabla^\perp u_p$.
Therefore,
\[
\Delta_{p'} v_p = 0.
\]
Thus we have the same sort of equation, but we can now consider the
limit $p' \to 1$. The duality between these two problems has been
exploited for different purposes before \cite{Aronsson-Lindqvist:88, Lindqvist:88}, but the consequences for the limit
behaviour have never been studied in detail, perhaps because swapping
$p \to \infty$ for $p' \to 1$ does not seem helpful superficially.
Here, however, is where the inverse mean curvature flow and its
$p'$-approximation come into play.

Set $w_p = \frac{1}{1 - p} \log v_p$. Then we compute
\[
\Delta_{p'} w_p = |\nabla w_p|^{p'}.
\]
Equation \eqref{eqn:imcf} arises as the formal limit as $p' \to 1$.
This connection between $p'$-harmonic functions and the inverse
mean curvature flow has been used before to construct weak
solutions of the latter \cite{Moser:07.1, Moser:08.1, Kotschwar-Ni:09, Moser:15.2}.
In the context of $\infty$-harmonic functions, the beauty in this
transformation is that it removes some of the degenerate behaviour
that arises for $v_p$ in the limit.

Next we use some tools developed for the inverse
mean curvature flow \cite{Moser:07.1, Moser:15.2} to show
that we have at least a sequence
$p_k \to \infty$ such that $w_{p_k}$ converges weakly in
$W_\loc^{1, q}(\Omega)$, for any $q < \infty$, to a weak solution $w$ of \eqref{eqn:imcf}.
Then we can reverse the above transformations to see what this means for
$u_\infty$. We compute
\[
|\nabla w_p|^{p' - 2}\nabla w_p = -(p' - 1)^{p' - 1} \omega e^{w_p} \nabla^\perp u_p.
\]
The left-hand side, at least if restricted to a certain subsequence,
will converge weakly in $L_\loc^q(\Omega; R^2)$ for every $q < \infty$
to a vector field $F$ satisfying the conditions from Definition \ref{def:imcf}.
The right-hand side converges to $-\omega e^w \nabla^\perp u_\infty$.
At almost every point $x \in \Omega$ such that $\nabla w(x) \neq 0$,
we conclude that
\[
\frac{\nabla w(x)}{|\nabla w(x)|} = -\omega e^{w(x)} \nabla^\perp u_\infty(x),
\]
and at such a point we therefore recover the relationship
$w(x) = -\log |\nabla u_\infty(x)|$ and also equations
\eqref{eqn:gradient} and \eqref{eqn:div}.

But it is possible that $\nabla w$ vanishes, and this is indeed
expected for situations such as when $\nabla u_\infty$ is constant.
In this case, we need much better information about the functions
$w_p$, and this is the most intricate part of the proof. We do not
go into the details here, but because of the technical difficulties
arising when $\nabla w = 0$, we will first consider
a small neighbourhood of a given point where
$\nabla u_\infty$ is nearly constant. As a consequence, we need
to show at the end of the proof
that a local weak solution of the inverse mean curvature
flow gives rise to a global weak solution. This is the main reason
why we favour Definition \ref{def:imcf} over the definition
of Huisken and Ilmanen \cite{Huisken-Ilmanen:01}. At least in the presence
of equations \eqref{eqn:gradient} and \eqref{eqn:div}, this step turns out to be
quite straightforward.

This strategy resembles some arguments that have been used for several
higher order variational problems related to the Aronsson equation
\cite{Moser-Schwetlick:12, Sakellaris:17, Katzourakis-Moser:19, Katzourakis-Parini:17, Moser:19}. These papers study minimisers of certain functionals involving
the $L^\infty$-norm. They rely on the idea of approximating the
$L^\infty$-norm by the $L^p$-norm for $p < \infty$, studying minimisers
of the resulting functionals, and reformulating the Euler-Lagrange
equation in a way that removes the expected degeneracy
in the limit $p \to \infty$, so that conclusions about the original
problem can be drawn. It is typically quite easy to find bounds
for the relevant quantities in the appropriate spaces in this context,
but it is necessary and difficult to show that they stay away from $0$.

It may seem that the above observations are specific to two-dimensional
domains, but they conceivably have a higher-dimensional
generalisation---not for the Aronsson equation \eqref{eqn:infty-harmonic},
but for an analogous problem involving differential forms.
Indeed, the relationship between $p'$-harmonic functions and
the $p'$-approximation of equation \eqref{eqn:imcf} exists for any dimension.
If $d$ denotes the exterior derivative and $d^*$ its formal $L^2$-adjoint,
then we may write the equation $\Delta_{p'} v_p = 0$ in the form
$d^*(|dv_p|^{p' - 2} dv_p) = 0$. Assuming that this is satisfied in
a star-shaped domain $\Omega \subseteq \R^n$, it implies that
$|dv_p|^{p' - 2} dv_p = d^*u_p$ for some $2$-form $u_p$ on $\Omega$. Moreover,
$u_p$ will satisfy $d(|d^*u_p|^{p - 2} d^*u_p) = 0$, which is the
Euler-Lagrange equation for the functional
\[
E_p(u) = \left(\frac{1}{p} \int_\Omega |d^* u|^p \, dx\right)^{1/p}.
\]
This suggests that we study the problem of minimising
$\|d^* u\|_{L^\infty(\Omega)}$ if we wish to find a connection to
the inverse mean curvature flow.
(If $n = 3$, we may alternatively minimise $\|\curl u\|_{L^\infty(\Omega)}$
for vector fields $u \colon \Omega \to \R^3$.)
Indeed, formal calculations
analogous to Aronsson's \cite{Aronsson:67} lead to the equation
\begin{equation} \label{eqn:2-form}
d|d^* u_\infty|^2 \wedge d^* u_\infty = 0.
\end{equation}
(For $n = 3$, we alternatively have the equation
$\nabla |\curl u_\infty|^2 \times \curl u_\infty = 0$.)
But almost nothing is known about this equation; indeed, even
the vector-valued optimal Lipschitz extension problem and
the Aronsson
equation for vector-valued functions $u \colon \Omega \to \R^N$
with $N \ge 2$ are poorly understood despite some existing work
on the former by Sheffield and Smart \cite{Sheffield-Smart:12}
and a series of papers on the latter by Katzourakis
\cite{Katzourakis:12, Katzourakis:13, Katzourakis:14.2, Katzourakis:14.1, Katzourakis:15.2, Katzourakis:17}. In particular, several of the tools for the proof
of Theorem \ref{thm:imcf} are missing in higher dimensions,
and we have no results here apart from the following calculations for
$C^2$-solutions, which are completely analogous to the above
calculations for $n = 2$.

Suppose that $u_\infty$ is a $2$-form with coefficients
in $C^2(\Omega)$. If $u_\infty$ solves \eqref{eqn:2-form} and
satisfies $d^*u_\infty \neq 0$ and $d|d^*u_\infty| \neq 0$ in $\Omega$,
then we conclude that
\[
\frac{d|d^* u_\infty|}{|d|d^* u_\infty||} = \pm\frac{d^*u_\infty}{|d^* u_\infty|}.
\]
Define $w = -\log|d^* u_\infty|$. Then
\[
-d^*\left(\frac{dw}{|dw|}\right) = \pm d^*\left(\frac{d^*u_\infty}{|d^* u_\infty|}\right) = \pm \frac{d|d^* u_\infty| \cdot d^* u_\infty}{|d^* u_\infty|^2} = \frac{|d|d^* u_\infty||}{|d^* u_\infty|} = |dw|.
\]
As the operator $-d^*$ for $1$-forms can be identified with the
divergence for vector fields, this means that $w$ solves equation \eqref{eqn:imcf}. Of course it is no longer appropriate to speak of streamlines
here. Their higher-dimensional counterparts are the hypersurfaces
characterised by the condition that their tangent vectors
$X$ satisfy $d^* u_\infty(X) = 0$, and using \eqref{eqn:2-form} we can check that
they coincide with the level sets of $|\nabla u_\infty|$.

The next few sections are devoted to the proofs of Theorem \ref{thm:imcf}
and Theorem \ref{thm:regularity}.
Then, in Section \ref{sct:imcf}, we prove that
weak solutions of \eqref{eqn:imcf} in the sense of Definition \ref{def:imcf}
are also weak solutions in the sense of Huisken and Ilmanen \cite{Huisken-Ilmanen:01}.
This final section is not essential for the understanding of the main
theorems, but it provides a connection with a larger body of literature
on the inverse mean curvature flow.

\section{Reduction to a local result}

As discussed in the introduction, we first consider small neighbourhoods
of a given point $x_0 \in \Omega$ where $\nabla u_\infty$ is nearly constant.
We may then rescale these neighbourhoods and thereby renormalise
$\nabla u_\infty(x_0)$
to unit size, using the following observation:
if $u_\infty$ is a given $\infty$-harmonic
function, then for any $a \in \R$, $r > 0$, and
$R \in \mathrm{O}(2)$, the rescaled function
$\tilde{u}_\infty(x) = au_\infty(rRx + x_0)$ is also
$\infty$-harmonic. If $a \det R > 0$, then the transformation preserves
the orientation, and if $a\det R < 0$, it reverses the orientation of
$u_\infty$. (We do not consider the case $a = 0$.)

The following result should be thought of as a statement about $u_\infty$
after such a rescaling, chosen such that $\nabla u_\infty(x_0)$ becomes
the second standard basis vector and the orientation becomes negative.
Here and throughout the rest of the paper, we use
the notation $(e_1, e_2)$ for the standard basis of $\R^2$, and we also write
$Q_r = (-r, r)^2$ for $r > 0$.

\begin{proposition} \label{prp:local}
There exists $\delta > 0$ with the following property.
Suppose that
$u_\infty \in C^1(\overline{Q_1})$ is $\infty$-harmonic with
$|\nabla u_\infty - e_2| \le \delta$ in $\overline{Q_1}$.
For $t \in [-\frac{1}{2}, \frac{1}{2}]$, let
$L_t = \set{x \in Q_1}{u_\infty(x) = u_\infty(0, t)}$, and suppose that
the numbers $m_t \in [-1, 1]$ satisfy the following condition: for
all $x, y \in L_t$, if $x_1 \le y_1 \le m_t$, then $|\nabla u_\infty(x)| \le |\nabla u_\infty(y)|$.
Let
\[
M = \bigcup_{t \in [-\frac{1}{2}, \frac{1}{2}]} \set{x \in L_t}{x_1 \le m_t}.
\]
Then there exists $w \in \bigcap_{q < \infty} W^{1,q}(Q_{1/4})$
such that
\begin{enumerate}
\item \label{itm:1} $w \le -\log |\nabla u_\infty|$ and
$e^w \nabla^\perp u_\infty \cdot \nabla w = |\nabla w|$ almost everywhere in $Q_{1/4}$,
\item \label{itm:2} $w = -\log |\nabla u_\infty|$ in $Q_{1/4} \cap M$, and
\item \label{itm:3} the equation
\[
\div (e^w \nabla^\perp u_\infty)= |\nabla w|
\]
holds weakly in $Q_{1/4}$.
\end{enumerate}
\end{proposition}

We give the proof of this result in Section \ref{sct:proof_local}
after some auxiliary results in Section \ref{sct:p-harmonic}.
But first, we show how Theorem \ref{thm:imcf} and Theorem \ref{thm:regularity}
follow from Proposition \ref{prp:local}.

\begin{proof}[Proof of Theorem \ref{thm:imcf}]
For any $x_0 \in \Omega$, we can choose $a \in \R$, $r > 0$, and
$R \in \mathrm{SO}(2)$ such that Proposition \ref{prp:local} applies
to $\tilde{u}(x) = au_\infty(rRx + x_0)$. Since $|\nabla \tilde{u}_\infty|$
is monotone along the level sets of $u_\infty$, we may choose
$m_t = 1$ for every $t \in [-\frac{1}{2}, \frac{1}{2}]$. Hence we obtain
a function $\tilde{w} \in \bigcap_{q < \infty} W^{1,q}(Q_{1/4})$ satisfying \ref{itm:1}--\ref{itm:3} in $Q_{1/4}$. In
particular
$\tilde{w} = -\log |\nabla \tilde{u}_\infty|$, and it follows that
$|\nabla \tilde{u}_\infty| \in W^{1, q}(Q_{1/4})$ for
every $q< \infty$.

From the pointwise equations \ref{itm:1} and \ref{itm:2}, we obtain
\[
\frac{\nabla^\perp \tilde{u}_\infty}{|\nabla \tilde{u}_\infty|} = -\frac{\nabla |\nabla \tilde{u}_\infty|}{|\nabla |\nabla \tilde{u}_\infty||}
\]
at almost every point where $\nabla |\nabla \tilde{u}_\infty| \neq 0$. This amounts to equation
\eqref{eqn:gradient} for $\tilde{u}_\infty$, which is trivially satisfied where the gradient
vanishes.
The combination of \ref{itm:2} and
\ref{itm:3} gives \eqref{eqn:div} for $\tilde{u}_\infty$. In terms of $u_\infty$,
this means that there exists a
neighbourhood $U$ of $x_0$ such that $|\nabla u_\infty| \in W^{1, q}(U)$
for every $q < \infty$ and \eqref{eqn:gradient} holds almost
everywhere in $U$, while
\eqref{eqn:div} holds weakly in $U$.

It follows that $|\nabla u_\infty| \in W_\loc^{1, q}(\Omega)$ and the
two equations hold in $\Omega$.
Let
\[
F = -\omega \frac{\nabla^\perp u_\infty}{|\nabla u_\infty|}.
\]
For the
function $w = - \log|\nabla u_\infty|$, we then compute
\[
\nabla w \cdot F = |\nabla w|
\]
because of \eqref{eqn:gradient},
and
\[
\div F = |\nabla w|
\]
because of \eqref{eqn:div}.
Thus $w$ is a weak solution of \eqref{eqn:imcf}.
\end{proof}

\begin{proof}[Proof of Theorem \ref{thm:regularity}]
For points in $G$, the arguments in the proof
of Theorem \ref{thm:imcf} apply. For $x_0 \in \Gamma$, we can still argue
similarly. Here we can still choose $a \in \R$, $r > 0$, and
$R \in \mathrm{SO}(2)$ such that the function
$\tilde{u}(x) = au_\infty(rRx + x_0)$ satisfies
$|\nabla \tilde{u}_\infty - e_2| < \delta$ in $\overline{Q_1}$
and such that
\[
\frac{1}{r}R^{-1}(G - x_0) \cap Q_1 = \set{x \in Q_1}{x_1 < f(x_2)}
\]
for some Lipschitz function $f \colon [-1, 1] \to \R$ with $f(0) = 0$,
the Lipschitz constant of which is independent of the rescaling.
We define $L_t$ as in Proposition \ref{prp:local} for $\tilde{u}_\infty$.
If $\delta$ is sufficiently small, then each $L_t$
will intersect the graph of $f$ exactly once for every $t \in [-\frac{1}{2}, \frac{1}{2}]$.
We choose $m_t$ such that the unique point $x \in L_t$ with $x_1 = m_t$
is this intersection point.

The orientation in Theorem \ref{thm:regularity} is such that
$|\nabla u_\infty|$ is non-decreasing if we approach $\Gamma$ from
inside $G$. In terms of $\tilde{u}_\infty$, this means that
$|\nabla \tilde{u}_\infty|$ is non-decreasing when we travel along
$L_t$ from left to right up to $m_t$. Then the hypothesis of
Proposition \ref{prp:local} is satisfied, so we infer $|\nabla \tilde{u}_\infty| \in W^{1, q}(M \cap Q_{1/4})$ for every $q < \infty$. Since $M \cap Q_{1/4}$
corresponds to a neighbourhood of $x_0$ in $G \cup \Gamma$, a standard covering argument now implies the desired statement.
\end{proof}

\section{Some estimates for $p$-harmonic functions} \label{sct:p-harmonic}

In this section we consider solutions of the equation $\Delta_p u = 0$.
We first prove an $L^\infty$-estimate for the gradient away from the boundary.
Such estimates are well known for fixed values of $p$, but we need to
know what happens when $p \to \infty$.
We use the well-known fact that $|\nabla u|^p$ is a subsolution to an elliptic
equation, and we derive an inequality with
the Moser iteration method. Here it is convenient to consider
higher-dimensional domains $\Omega \subseteq \R^n$ for $n \ge 3$,
as the proof uses the Sobolev embedding
$W_0^{1, 2}(\Omega) \subset L^{\frac{2n}{n - 2}}(\Omega)$,
which does not exist in this form for $n = 2$. Of course we still obtain
an estimate for $\Omega \subseteq \R^2$ by extending a $p$-harmonic function
$u \colon \Omega \to \R$ to $\Omega \times (0, 1)$ by
$\tilde{u}(x_1, x_2, x_3) = u(x_1, x_2)$.

\begin{lemma} \label{lem:L^infty}
For any $n \ge 3$ there exists a constant $C > 0$ with
the following property. Suppose that $\Omega \subseteq \R^n$ is an open set.
For $\varrho > 0$, let $\Omega_\varrho = \set{x \in \Omega}{\dist(x, \partial \Omega) > \varrho}$. Let $p \ge 2$.
Then for any $p$-harmonic function $u \in W^{1, p}(\Omega)$,
\[
\|\nabla u\|_{L^\infty(\Omega_\varrho)} \le \left(\frac{C p^{\frac{n}{2} + 1}}{\varrho^n}\right)^{\frac{1}{p}} \|\nabla u\|_{L^p(\Omega)}.
\]
\end{lemma}

\begin{proof}
First we remark that it suffices to consider solutions
of the equation
\begin{equation} \label{eqn:approx-p-harmonic}
\div\left(\bigl(|\nabla u|^2 + \epsilon^2\bigr)^{\frac{p}{2} - 1} \nabla u\right) = 0
\end{equation}
in $\Omega$ for $\epsilon > 0$ and to prove the inequality
\begin{equation} \label{eqn:pointwise_epsilon}
\left\|\sqrt{|\nabla u|^2 + \epsilon^2}\right\|_{L^\infty(\Omega_\varrho)} \le \left(\frac{C p^{\frac{n}{2} + 1}}{\varrho^n}\right)^{\frac{1}{p}} \left\|\sqrt{|\nabla u|^2 + \epsilon^2}\right\|_{L^p(\Omega)}.
\end{equation}
The statement of the lemma then follows by approximation arguments as discussed, e.g., by Lewis \cite{Lewis:83}.
Thus we assume that $u$ satisfies \eqref{eqn:approx-p-harmonic}.
Well-known regularity theory \cite[Chapter 4]{Ladyzhenskaya-Uraltseva:68}
then shows that $u$ is smooth.

Let $g = \sqrt{|\nabla u|^2 + \epsilon^2}$. We differentiate
with respect to $x_i$ in \eqref{eqn:approx-p-harmonic}, which gives
\[
\div\left(g^{p - 2} \nabla \dd{u}{x_i} + (p - 2) g^{p - 4} \left(\nabla u \cdot \nabla \dd{u}{x_i}\right) \nabla u\right) = 0.
\]
Multiply with $\dd{u}{x_i}$. This gives
\begin{multline} \label{eqn:differentiated}
\div\left(g^{p - 2} \dd{u}{x_i} \nabla\dd{u}{x_i} + (p - 2) g^{p - 4} \dd{u}{x_i} \left(\nabla u \cdot \nabla \dd{u}{x_i}\right) \nabla u\right) \\
= g^{p - 2} \left|\nabla \dd{u}{x_i}\right|^2 + (p - 2) g^{p - 4} \left(\nabla u \cdot \nabla \dd{u}{x_i}\right)^2.
\end{multline}
For any $\eta \in C_0^\infty(\Omega)$, we then obtain
\begin{multline*}
\int_\Omega \eta^2 \left(g^{p - 2} \left|\nabla \dd{u}{x_i}\right|^2 + (p - 2) g^{p - 4} \left(\nabla u \cdot \nabla \dd{u}{x_i}\right)^2\right) \, dx \\
\begin{aligned}
& = - 2\int_\Omega \eta \nabla \eta \cdot \left(g^{p - 2} \dd{u}{x_i} \nabla  \dd{u}{x_i} + (p - 2) g^{p - 4} \dd{u}{x_i} \left(\nabla u \cdot \nabla \dd{u}{x_i}\right) \nabla u\right) \, dx \\
& \le \frac{1}{2} \int_\Omega \eta^2 \left(g^{p - 2} \left|\nabla \dd{u}{x_i}\right|^2 + (p - 2) g^{p - 4} \left(\nabla u \cdot \nabla \dd{u}{x_i}\right)^2\right) \, dx \\
& \quad + 2(p - 1) \int_\Omega |\nabla \eta|^2 g^{p - 2} \left(\dd{u}{x_i}\right)^2 \, dx.
\end{aligned}
\end{multline*}
Hence
\begin{multline*}
\int_\Omega \eta^2 \left(g^{p - 2} \left|\nabla \dd{u}{x_i}\right|^2 + (p - 2) g^{p - 4} \left(\nabla u \cdot \nabla \dd{u}{x_i}\right)^2\right) \, dx\\
\le  4(p - 1) \int_\Omega |\nabla \eta|^2 g^{p - 2} \left(\dd{u}{x_i}\right)^2 \, dx.
\end{multline*}
It follows that
\begin{equation} \label{eqn:initial_inequality}
\int_\Omega \eta^2 \left(\dd{}{x_i}g^{p/2}\right)^2 \, dx \le p^2 \int_\Omega |\nabla \eta|^2 g^{p - 2} \left(\dd{u}{x_i}\right)^2 \, dx.
\end{equation}

Another consequence of \eqref{eqn:differentiated} is the following: define
\[
f = \frac{g^p}{p}.
\]
Let $I$ denote the identity $(n \times n)$-matrix and consider
\[
A = I + (p - 2) g^{-2} \nabla u \otimes \nabla u.
\]
Then we compute
\[
\div (A\nabla f) \ge \frac{p - 1}{p} \frac{|\nabla f|^2}{f}.
\]
For $\eta \in C_0^\infty(\Omega)$ and $q > -1$, we therefore obtain the inequality
\[
\begin{split}
\lefteqn{\int_\Omega \eta^2 f^q \nabla f \cdot A\nabla f \, dx} \qquad \\
& \le -\frac{2}{q + 1} \int_\Omega \eta f^{q + 1} \nabla \eta \cdot A\nabla f \, dx - \frac{p - 1}{p(q + 1)} \int_\Omega \eta^2 f^q |\nabla f|^2\, dx \\
& \le \frac{1}{2} \int_\Omega \eta^2 f^q \nabla f \cdot A\nabla f \, dx + \frac{2}{(q + 1)^2} \int_\Omega f^{q + 2} \nabla \eta \cdot A\nabla \eta \, dx \\
& \quad - \frac{p - 1}{p(q + 1)} \int_\Omega \eta^2 f^q |\nabla f|^2\, dx.
\end{split}
\]
Hence
\begin{multline*}
\int_\Omega \eta^2 f^q \nabla f \cdot A\nabla f \, dx + \frac{2(p - 1)}{p(q + 1)} \int_\Omega \eta^2 f^q |\nabla f|^2\, dx \\
\le \frac{4}{(q + 1)^2} \int_\Omega f^{q + 2} \nabla \eta \cdot A\nabla \eta \, dx.
\end{multline*}
In particular,
\begin{equation} \label{eqn:iterated_inequality}
\left(1 + \frac{2(p - 1)}{p(q + 1)}\right) \int_\Omega \eta^2 f^q |\nabla f|^2\, dx
\le \frac{4(p - 1)}{(q + 1)^2} \int_\Omega f^{q + 2} |\nabla \eta|^2 \, dx.
\end{equation}

Now choose $s > \frac{1}{2}$. By the Sobolev inequality, we have the estimate
\[
\begin{split}
\left(\int_\Omega \eta^{\frac{2n}{n - 2}} f^{\frac{2ns}{n - 2}} \, dx\right)^{\frac{n - 2}{2n}} & \le C_1 \left(\int_\Omega |\nabla  (\eta f^s)|^2 \, dx\right)^{\frac{1}{2}} \\
& \le C_1\left(\int_\Omega f^{2s} |\nabla \eta|^2 \, dx\right)^{\frac{1}{2}} + C_1s \left(\int_\Omega \eta^2 f^{2s - 2} |\nabla f|^2 \, dx\right)^{\frac{1}{2}}
\end{split}
\]
for some constant $C_1 = C_1(n)$.
Using \eqref{eqn:iterated_inequality} for $q = 2s - 2$, we obtain
\[
\begin{split}
\int_\Omega \eta^2 f^{2s - 2} |\nabla f|^2 \, dx & \le \frac{4(p - 1)}{(2s - 1) \left(2s -1 + \frac{2(p - 1)}{p}\right)} \int_\Omega f^{2s} |\nabla \eta|^2 \, dx \\
& \le \frac{2(p - 1)}{s(2s - 1)} \int_\Omega f^{2s} |\nabla \eta|^2 \, dx,
\end{split}
\]
where we have used the assumption that $p \ge 2$. Hence
\[
\left(\int_\Omega \eta^{\frac{2n}{n - 2}} f^{\frac{2ns}{n - 2}} \, dx\right)^{\frac{n - 2}{2n}} \le C_1\left(1 + \sqrt{\frac{2s(p - 1)}{2s - 1}}\right) \left(\int_\Omega f^{2s} |\nabla \eta|^2 \, dx\right)^{\frac{1}{2}}.
\]

An appropriate choice of $\eta$ in this inequality gives
\[
\left(\int_{\Omega_R} f^{\frac{2ns}{n - 2}} \, dx\right)^{\frac{n - 2}{2n}} \le \frac{C_1}{R - r} \left(1 + \sqrt{\frac{2s(p - 1)}{2s - 1}}\right) \left(\int_{\Omega_r} f^{2s} \, dx\right)^{\frac{1}{2}}
\]
whenever $0 < r < R$. In particular, for any $q_0 > 1$, there exists
a constant $C_2 = C_2(n, q_0)$ such that
\[
\left(\int_{\Omega_R} f^{\frac{2ns}{n - 2}} \, dx\right)^{\frac{n - 2}{2n}} \le \frac{C_2 \sqrt{p}}{R - r} \left(\int_{\Omega_r} f^{2s} \, dx\right)^{\frac{1}{2}}
\]
for all $s \ge q_0/2$.

Define
\[
J(q, r) = \left(\int_{\Omega_r} f^q \, dx\right)^{\frac{1}{q}}
\]
for $q \ge q_0$ and $r > 0$. Set $\theta = \frac{n}{n - 2}$. Then the preceding inequality can be
written as
\[
J(\theta q, R) \le \left(\frac{C_2 \sqrt{p}}{R - r}\right)^{\frac{2}{q}} J(q, r)
\]
for $q \ge q_0$. Set $r_k = \varrho_0 + (1 - 2^{-k}) (\varrho - \varrho_0)$ for some $\varrho_0, \varrho > 0$
with $\varrho > \varrho_0$ and set $q_k = \theta^k q_0$.
Then
\[
J(q_{k + 1}, r_{k + 1}) \le \left(\frac{2^{k + 1} C_2 \sqrt{p}}{\varrho - \varrho_0}\right)^{\frac{2}{q_k}} J(q_k, r_k).
\]
Iterating this, we obtain
\[
\|f\|_{L^\infty(\Omega_\varrho)} \le C_3 \|f\|_{L^{q_0}(\Omega_{\varrho_0})},
\]
where
\[
C_3 = \prod_{k = 0}^\infty \left(\frac{2^{k + 1} C_2 \sqrt{p}}{\varrho - \varrho_0}\right)^{\frac{2}{q_k}}.
\]
We can rewrite $C_3$ as follows: let
\[
a = \sum_{k = 0}^\infty \frac{1}{q_k} = \frac{1}{q_0} \sum_{k = 0}^\infty \theta^{-k} = \frac{\theta}{q_0(\theta - 1)} = \frac{n}{2q_0}
\]
and
\[
b = \sum_{k = 0}^\infty \frac{k + 1}{q_k} = \frac{1}{q_0} \sum_{k = 0}^\infty (k + 1) \theta^{-k} = \frac{\theta^2}{q_0 (\theta - 1)^2} = \frac{n^2}{4q_0}.
\]
Then
\[
C_3 = 4^b \left(\frac{C_2^2 p}{(\varrho - \varrho_0)^2}\right)^a = \left(\frac{2^{n^2/2} C_2^n p^{n/2}}{(\varrho - \varrho_0)^n}\right)^{\frac{1}{q_0}}.
\]

We can finally do one more step using \eqref{eqn:initial_inequality}. Using the Sobolev
inequality the same way as before, we obtain
\[
\begin{split}
\left(\int_\Omega \left(\eta g^{p/2}\right)^{\frac{2n}{n - 2}} \, dx\right)^{\frac{n - 2}{2n}} & \le C_1\left(\int_\Omega |\nabla \eta|^2 g^p \, dx\right)^{\frac{1}{2}} + C_1\left(\int_\Omega \eta^2 \left|\nabla g^{p/2}\right|^2 \, dx\right)^{\frac{1}{2}} \\
& \le C_1 (p + 1) \left(\int_\Omega |\nabla \eta|^2 g^p \, dx\right)^{\frac{1}{2}}.
\end{split}
\]
Hence if we use suitable functions $\eta$ again, then in terms of $f$,
we derive the inequality
\[
\left(\int_{\Omega_{\varrho_0}} f^{\frac{n}{n - 2}} \, dx\right)^{\frac{n - 2}{n}} \le \frac{C_1^2}{\varrho_0^2} (p + 1)^2 \int_\Omega f \, dx.
\]
Choose $\varrho_0 = \varrho/2$ and $q_0 = \frac{n}{n - 2}$. Then it follows that
\[
\|f\|_{L^\infty(\Omega_\varrho)} \le \frac{C_4 p^{(n + 2)/2}}{\varrho^n} \|f\|_{L^1(\Omega)},
\]
for some constant $C_4 = C_4(n)$. In terms of $u$, this is
inequality \eqref{eqn:pointwise_epsilon}.
\end{proof}

We go back to domains in the plane. The following result gives
an estimate from below for the partial derivative $\dd{u}{x_2}$
of a $p$-harmonic function in a suitable domain with suitable boundary data.

\begin{lemma} \label{lem:partial_derivative}
Let $f, g \colon [-1, 1] \to \R$ be two Lipschitz functions with
$f < g$ and $U = \set{x \in (-1, 1) \times \R}{f(x_1) < x_2 < g(x_1)}$.
Let $\phi \colon \overline{U} \to \R$ be a Lipschitz function
such that $\dd{\phi}{x_2} \ge 1$ for almost all $x \in U$
and such that there are two numbers $a, b \in \R$ with $\phi(x) = a + x_2$
when $x_2 = f(x_1)$ and $\phi(x) = b + x_2$ when $x_2 = g(x_1)$.
Then the solution $u \colon \overline{U} \to \R$ of the boundary
value problem
\begin{alignat*}{2}
\Delta_p u & = 0 & \quad & \text{in $U$}, \\
u & = \phi && \text{on $\partial U$},
\end{alignat*}
satisfies $\dd{u}{x_2} \ge 1$ in $U$.
\end{lemma}

\begin{proof}
Let $c = \max\{\|f\|_{L^\infty(-1, 1)}, \|g\|_{L^\infty(-1, 1)}\} + 2$
and write $U_0 = (-1, 1) \times (-c, c)$.
Extend $\phi$ to $U_0$ by $\phi(x) = a + x_2$ when $x_2 < f(x_1)$
and $\phi(x) = b + x_2$ when $x_2 > g(x_1)$.
Choose a sequence of functions
$\phi_k \in C^\infty(\overline{U_0})$, for $k \in \N$,
such that $\phi_k \to \phi$
in $W^{1, p}(U_0)$ as $k \to \infty$ and such that
\begin{itemize}
\item $\phi_k(x) = a + x_2$ when $x_2 < f(x_1) - \frac{1}{k}$,
\item $\phi_k(x) = b + x_2$ when $x_2 > g(x_1) + \frac{1}{k}$, and
\item $\dd{\phi_k}{x_2} \ge 1$.
\end{itemize}
Now choose a sequence
of domains $U_k \subseteq U_0$ with smooth boundaries,
such that each $U_k$ is of the form
$U_k = \set{x \in (-1, 1) \times \R}{f_k(x_1) < x_2 < g_k(x_1)}$
for some smooth functions $f_k, g_k \colon (-1, 1) \to \R$ with
\[
f(x_1) -\frac{2}{k} < f_k(x_1) < f(x_1) - \frac{1}{k} < g(x_1) + \frac{1}{k} < g_k(x_1) < g(x_1) + \frac{2}{k}
\]
for $-1 < x_1 < 1$.

Let $u_k \colon U_k \to \R$ be the solution of
\begin{alignat*}{2}
\div\left(\bigl(|\nabla u_k|^2 + k^{-2}\bigr)^{p/2 - 1} \nabla u_k\right) & = 0 & \quad & \text{in $U_k$}, \\
u_k & = \phi_k && \text{on $\partial U_k$}.
\end{alignat*}
Extend $u_k$ to $U_0$ by $u_k(x) = a + x_2$ when
$x_2 < f_k(x_1)$ and $u_k(x) = b + x_2$ when $x_2 > g_k(x_1)$.
Then the sequence $(u_k)_{k \in \N}$ is clearly bounded in
$W^{1, p}(U_0)$, and we may assume that it converges
weakly in this space to a limit $\tilde{u} \in W^{1, p}(U_0)$.
We claim that $\tilde{u} = u$ in $U$.

In order to prove this, note that
\[
\int_{U_k} \left(|\nabla u_k|^2 + k^{-2}\right)^{p/2} \, dx \le \int_{U_k} \left(|\nabla (\phi_k + w)|^2 + k^{-2}\right)^{p/2} \, dx
\]
for any $w \in W_0^{1, p}(U)$, because $u_k$ minimises this
quantity for its boundary data.
Letting $k \to \infty$, we find that
\[
\begin{split}
\int_U |\nabla \tilde{u}|^p \, dx & \le \liminf_{k \to \infty} \int_{U_k} \left(|\nabla u_k|^2 + k^{-2}\right)^{p/2} \, dx \\
& \le \liminf_{k \to \infty} \int_{U_k} \left(|\nabla (\phi_k + w)|^2 + k^{-2}\right)^{p/2} \, dx.
\end{split}
\]
Moreover,
\begin{multline*}
\int_{U_k} \left(|\nabla (\phi_k + w)|^2 + k^{-2}\right)^{p/2} \, dx \\
\begin{aligned}
& = \int_{U_0} \left(|\nabla (\phi_k + w)|^2 + k^{-2}\right)^{p/2} \, dx - \int_{U_0 \setminus U_k} \left(|\nabla \phi|^2 + k^{-2}\right)^{p/2} \, dx \\
& \to \int_{U_0} |\nabla (\phi + w)|^p \, dx - \int_{U_0 \setminus U} |\nabla \phi|^p \, dx \\
& = \int_U |\nabla (\phi + w)|^p \, dx
\end{aligned}
\end{multline*}
as $k \to \infty$. Therefore,
\[
\int_U |\nabla \tilde{u}|^p \, dx \le \int_U |\nabla (\phi + w)|^p \, dx
\]
for any $w \in W_0^{1, p}(U)$. Furthermore, it is clear that
$\tilde{u} = \phi$ on $\partial U$. Since the functional
\[
v \mapsto \int_U |\nabla v|^p \, dx
\]
has a unique minimiser in $\phi + W_0^{1, p}(U)$, which is $u$, we
conclude that $\tilde{u} = u$.

The regularity theory of Lieberman \cite{Lieberman:88} shows that
$u_k \in C^\infty(\overline{U_k})$.
Moreover, the function $\dd{u_k}{x_2}$ satisfies the equation
\[
\div\left(A_k \nabla \dd{u_k}{x_2}\right) = 0,
\]
where
\[
A_k = \left(|\nabla u_k|^2 + k^{-2}\right)^{p/2 - 1} \left(I + (p - 2) \frac{\nabla u_k \otimes \nabla u_k}{|\nabla u_k|^2 + k^{-2}}\right)
\]
and $I$ is the identity matrix. This is a uniformly elliptic equation.

The comparison principle applies to $u_k$ and
implies that $a + x_2 \le u_k(x) \le b + x_2$ for
$x \in U_k$. Hence $\dd{u_k}{x_2} \ge 1$ on
\[
\set{x \in (-1, 1) \times \R}{x_2 = f_k(x_1) \text{ or } x_2 = g_k(x_1)}.
\]
On the rest of $\partial U_k$, the inequality is inherited
directly from $\phi_k$.
Applying the maximum principle to $\dd{u_k}{x_2}$, we prove that $\dd{u_k}{x_2} \ge 1$
in $U_k$. Then the desired inequality follows for $u$ as well.
\end{proof}

\section{Proof of Proposition \ref{prp:local}} \label{sct:proof_local}

This section contains the key arguments of this paper.

We first fix $\delta \in (0, \frac{1}{16})$ and also fix a constant $\sigma \in [\frac{1}{2}, 1 - 8\delta)$.
The values of both will be determined later. (We will require
that $\delta$ and $1 - \sigma$ are sufficiently small.)
We consider an $\infty$-harmonic function
$u_\infty \in C^1(\overline{Q_1})$ that satisfies the hypotheses of
Proposition \ref{prp:local}.
We wish to find a function $w \in \bigcap_{q < \infty} W^{1, q}(Q_{1/4})$
with the properties \ref{itm:1}--\ref{itm:3}. In particular, we wish
to show that $w$
is a weak solution of the inverse mean curvature flow.

The function $\tilde{u}(x) = u_\infty(x) - \sigma x_2$
satisfies
\[
|\nabla \tilde{u} - (1 - \sigma)e_2| \le \delta.
\]
Hence the level sets of $\tilde{u}$ are Lipschitz graphs with Lipschitz
constants bounded by
\[
\frac{\delta}{\sqrt{(1 - \sigma)^2 - \delta^2}}.
\]
Let $a = \tilde{u}(0, -\frac{3}{4})$ and $b = \tilde{u}(0, \frac{3}{4})$.
Set
\[
U = \set{x \in Q_1}{a < \tilde{u}(x) < b}.
\]
Since $\delta \le \frac{1}{8}(1 - \sigma)$ by the choice of $\sigma$,
it follows that
$[-1, 1] \times [-\frac{1}{2}, \frac{1}{2}] \subseteq
\overline{U} \subseteq [-1, 1] \times (-1, 1)$.

For $p \in [2, \infty)$, let $u_p \in W^{1, p}(U)$ be the unique
solutions of
\begin{alignat*}{2}
\Delta_p u_p & = 0 & \quad & \text{in $U$}, \\
u & = u_\infty && \text{on $\partial U$}.
\end{alignat*}
Lemma \ref{lem:partial_derivative} implies that
\begin{equation} \label{eqn:partial_derivative}
\dd{u_p}{x_2} \ge \sigma
\end{equation}
in $U$. Moreover,
\[
\int_U |\nabla u_p|^p \, dx \le \int_U |\nabla u_\infty|^p \, dx \le 4 (1 + \delta)^p,
\]
as $u_p$ minimises this quantity among all functions with the same
boundary data. Define
$U_p = \set{x \in U}{\dist(x, \partial U) > 2^{-\sqrt{p}}}$.
Then Lemma \ref{lem:L^infty} implies that
\[
\|\nabla u_p\|_{L^\infty(U_p)} \le \left(4C p^{\frac{5}{2}} 8^{\sqrt{p}}\right)^{\frac{1}{p}} (1 + \delta),
\]
where $C$ is the constant from Lemma \ref{lem:L^infty} for $n = 3$.
In particular, if $p$ is sufficiently large, then
\begin{equation} \label{eqn:L^infty}
\|\nabla u_p\|_{L^\infty(U_p)} \le 1 + 2\delta.
\end{equation}
Inequalities \eqref{eqn:partial_derivative} and \eqref{eqn:L^infty} then
imply that
\begin{equation} \label{eqn:angle_u}
\dd{u_p}{x_2}\ge \frac{\sigma |\nabla u_p|}{1 + 2\delta}.
\end{equation}
in $U_p$.

The results of Jensen \cite{Jensen:93} imply that $u_p \rightharpoonup u_\infty$
weakly in $W^{1, q}(U)$ for any $q < \infty$
and also uniformly in $U$. Moreover, we have the following variant of a
result by Lindgren and Lindqvist
\cite{Lindgren-Lindqvist:21}.

\begin{lemma} \label{lem:uniform_convergence}
For any precompact set $K \Subset U$, the convergence
$|\nabla u_p| \to |\nabla u_\infty|$ holds uniformly in $K$.
\end{lemma}

\begin{proof}
The arguments of Lindgren and Lindqvist \cite[Section 3]{Lindgren-Lindqvist:21}
(which depend to some degree on the ideas of Koch, Zhang, and
Zhou \cite{Koch-Zhang-Zhou:19.1} and also use an inequality of
Lebesgue \cite{Lebesgue:07})
can be used here. Although their paper deals with an annular domain and with
specific boundary conditions, their reasoning applies more generally to $p$-harmonic
functions in a domain $U \subseteq \R^2$ satisfying $|\nabla u_p| > 0$ in $U$ and
\[
\limsup_{p \to \infty} \|\nabla u_p\|_{L^\infty(K)} < \infty
\]
for any $K \Subset U$. In our case, the first property follows from
\eqref{eqn:partial_derivative} and the second from \eqref{eqn:L^infty}.
\end{proof}

\begin{remark}
Lemma \ref{lem:uniform_convergence} does not imply that
$\nabla u_p \to \nabla u_\infty$. Nevertheless, the ideas of Koch, Zhang,
and Zhou \cite{Koch-Zhang-Zhou:19.1}, as adapted by Lindgren and Lindqvist
\cite{Lindgren-Lindqvist:21}, do give convergence almost everywhere.
But we do not need this information here.
\end{remark}

For $p \in [2, \infty)$, let $p' \in (1, 2]$ denote the conjugate exponent
with $\frac{1}{p} + \frac{1}{p'} = 1$.
Since $U$ is simply connected and since
$\curl(|\nabla u_p|^{p - 2} \nabla^\perp u_p) = \Delta_p u_p = 0$,
there exists $v_p \in W^{1, p'}(U)$ satisfying
\[
\nabla v_p = -|\nabla u_p|^{p - 2} \nabla^\perp u_p.
\]
Then we compute $|\nabla v_p|^{p' - 2} \nabla v_p = -\nabla^\perp u_p$.
Therefore,
\[
\Delta_{p'} v_p = 0
\]
in $U$. Note that the level sets of $v_p$ are the streamlines of $u_p$.
Moreover, \eqref{eqn:angle_u} implies that
\begin{equation} \label{eqn:angle_v}
\dd{v_p}{x_1}\ge \frac{\sigma |\nabla v_p|}{1 + 2\delta}
\end{equation}
in $U_p$. That is, the angle between $\nabla v_p$ and $e_1$ is at most
\[
\arccos\left(\frac{\sigma}{1 + 2\delta}\right).
\]
If $\delta$ and $1 - \sigma$ are sufficiently small, then this means that
for every $\beta \in \R$ there exists $\theta \in [-1, 1]$ such that
\[
\textstyle (-1, \theta) \times [-\frac{1}{2}, \frac{1}{2}] \subseteq \set{x \in (-1, 1) \times [-\frac{1}{2}, \frac{1}{2}]}{v_p(x) < \beta} \subseteq (-1, \theta + \frac{1}{16}) \times [-\frac{1}{2}, \frac{1}{2}]
\]
for $p$ large enough.

We can further prove the following inequalities for $v_p$.
Here we use the notation
$a_+ = \max\{a, 0\}$ for $a \in \R$.
For $\xi = (\xi_1, \xi_2) \in \R^2$, we write $\xi^\perp = (-\xi_2, \xi_1)$.

\begin{lemma} \label{lem:pointwise}
Let $\lambda \colon [0, \ell] \to U$ be a $C^1$-curve with
$|\lambda'| \equiv 1$ and
$|\lambda' - e_1| \le \delta$. Let $x = \lambda(0)$ and $y = \lambda(\ell)$.
Then for any $\varrho > 0$,
\[
\limsup_{p \to \infty} \biggl(\inf_{B_\varrho(y) \cap U} v_p - \sup_{B_\varrho(x) \cap U} v_p\biggr)_+^{\frac{1}{p - 1}} \le \sup_{0 \le s \le \ell} |\nabla u_\infty(\lambda(s))|
\]
and
\[
\liminf_{p \to \infty} \biggl(\sup_{B_\varrho(y) \cap U} v_p - \inf_{B_\varrho(x) \cap U} v_p\biggr)^{\frac{1}{p - 1}} \ge \fint_0^\ell (\lambda'(s))^\perp \cdot \nabla u_\infty(\lambda(s)) \, ds.
\]
\end{lemma}

\begin{proof}
We may approximate $\lambda$ in the $C^1$-topology with smooth curves.
Therefore, we may assume without loss of generality that $\lambda \in C^\infty([0, \ell]; U)$.
Let $\epsilon > 0$.
Define $\Phi \colon [0, \ell] \times [- T, T] \to U$
by $\Phi(s, t) = \lambda(s) + t(\lambda'(s))^\perp$, where $T > 0$
is chosen so small that
\begin{equation} \label{eqn:curve_approx}
|\nabla u_\infty(\Phi(s, t)) - \nabla u_\infty(\lambda(s))| \le \epsilon
\end{equation}
for all $s \in [0, \ell]$ and $t \in [-T, T]$. Since
$\det D\Phi(s, t) = 1 + t\lambda'(s) \cdot (\lambda''(s))^\perp$, we may further
choose $T$ so small that $|\det D\Phi - 1| < \epsilon$ and
$|(\det D\Phi)^{-1} - 1| < \epsilon$ in
$[0, \ell] \times [-T, T]$. We write
$\Sigma = \Phi([0, \ell] \times [-T, T])$ and
$\Sigma_t = \Phi([0, \ell] \times \{t\})$ for $-T \le t \le T$,
and we set $X = \dd{\Phi}{t} \circ \Phi^{-1}$. (I.e., $X(\Phi(s, t)) = (\lambda'(s))^\perp$.)

Now note that $\Sigma \subseteq U_p$ for $p$ sufficiently large and $T$ sufficiently small.
Hence $X \cdot \nabla u_p \ge 0$ by \eqref{eqn:angle_u}
and $X^\perp \cdot \nabla v_p \le 0$ by \eqref{eqn:angle_v}.

We write $\Ha^1$ for the $1$-dimensional Hausdorff measure. Then
\begin{equation} \label{eqn:pntw1}
\begin{split}
\fint_{-T}^T \bigl(v_p(\Phi(\ell, t)) - v_p(\Phi(0, t))\bigr) \, dt
& = -\fint_{-T}^T \int_{\Sigma_t} X^\perp \cdot \nabla v_p \, d\Ha^1 \, dt \\
& = - \frac{1}{2T} \int_\Sigma X^\perp \cdot \nabla v_p \, dx \\
& = \frac{1}{2T} \int_\Sigma |\nabla u_p|^{p - 2} X \cdot \nabla u_p \, dx \\
& \le \frac{|\Sigma|}{2T} \sup_{\Sigma} |\nabla u_p|^{p - 1}.
\end{split}
\end{equation}
By \eqref{eqn:curve_approx} and Lemma \ref{lem:uniform_convergence},
if $p$ is sufficiently large,
then
\[
\sup_{\Sigma} |\nabla u_p| \le \sup_{0 \le s \le \ell} |\nabla u_\infty(\lambda(s))| + 2\epsilon.
\]
If $T < \varrho$, it follows that
\[
\biggl(\inf_{B_\varrho(y) \cap U} v_p - \sup_{B_\varrho(x) \cap U} v_p\biggr)_+^{\frac{1}{p - 1}} \le \left(\frac{|\Sigma|}{2T}\right)^{\frac{1}{p - 1}} \left(\sup_{0 \le s \le \ell} |\nabla u_\infty(\lambda(s))| + 2 \epsilon\right).
\]
Letting $p \to \infty$ and $\epsilon \to 0$, we obtain the first inequality.

We can also estimate
\begin{multline} \label{eqn:pntw2}
\int_0^\ell \bigl(u_p(\Phi(s, T)) - u_p(\Phi(s, -T))\bigr) \, ds \\
\begin{aligned}
& = \int_0^\ell \int_{-T}^T \dd{\Phi}{t}(s, t) \cdot \nabla u_p(\Phi(s, t)) \, dt \, ds \\
& = \int_\Sigma X \cdot \nabla u_p \, |\det D\Phi^{-1}| \, dx \\ 
& \le (1 + \epsilon) \int_\Sigma X \cdot \nabla u_p \, dx \\
& \le (1 + \epsilon) |\Sigma|^{\frac{p - 2}{p - 1}} \left(\int_\Sigma (X \cdot \nabla u_p)^{p - 1} \, dx\right)^{\frac{1}{p - 1}} \\
& \le (1 + \epsilon) |\Sigma|^{\frac{p - 2}{p - 1}} \left(\int_\Sigma |\nabla u_p|^{p - 2} X \cdot \nabla u_p \, dx\right)^{\frac{1}{p - 1}} \\
& = (1 + \epsilon) |\Sigma|^{\frac{p - 2}{p - 1}} \left(-\int_\Sigma X^\perp \cdot \nabla v_p \, dx\right)^{\frac{1}{p - 1}}.
\end{aligned}
\end{multline}
Since $u_p \to u_\infty$ uniformly, for $p$ sufficiently large we have
the inequality
\begin{multline*}
\int_0^\ell \bigl(u_p(\Phi(s, T)) - u_p(\Phi(s, -T))\bigr) \, ds \\
\begin{aligned}
& \ge (1 - \epsilon) \int_0^\ell \bigl(u_\infty(\Phi(s, T)) - u_\infty(\Phi(s, -T))\bigr) \, ds \\
& = \int_0^\ell \int_{-T}^T \dd{\Phi}{t}(s, t) \cdot \nabla u_\infty(\Phi(s, t)) \, dt \, ds.
\end{aligned}
\end{multline*}
Moreover, by \eqref{eqn:curve_approx} and because $(\lambda'(s))^\perp \cdot \nabla u_\infty(\lambda(s)) \ge \frac{1}{2}$ by the assumptions on $\lambda$,
we know that $\dd{\Phi}{t}(s, t) \cdot \nabla u_\infty(\Phi(s, t)) \ge (1 - 2\epsilon)(\lambda'(s))^\perp \cdot \nabla u_\infty(\lambda(s))$. Thus
\begin{equation} \label{eqn:pntw3}
\int_0^\ell \bigl(u_p(\Phi(s, T)) - u_p(\Phi(s, -T))\bigr) \, ds \ge 2\ell T (1 - 2\epsilon) \fint_0^\ell (\lambda'(s))^\perp \cdot \nabla u_\infty(\lambda(s)) \, ds.
\end{equation}

Recall that
\[
\fint_{-T}^T \bigl(v_p(\Phi(\ell, t)) - v_p(\Phi(0, t))\bigr) \, dt = -\frac{1}{2T} \int_\Sigma X^\perp \cdot \nabla v_p \, dx
\]
according to \eqref{eqn:pntw1}.
We also know that $|\Sigma| \le 2(1 + \epsilon)\ell T$.
Therefore, combining the above inequalities \eqref{eqn:pntw2} and
\eqref{eqn:pntw3}, we obtain
\begin{multline*}
\fint_{-T}^T \bigl(v_p(\Phi(\ell, t)) - v_p(\Phi(0, t))\bigr) \, dt \\
\ge\frac{(1 - 2\epsilon)^{p - 1}}{(1 + \epsilon)^{2p - 3}} \ell \left(\fint_0^\ell (\lambda'(s))^\perp \cdot \nabla u_\infty(\lambda(s)) \, ds\right)^{p - 1}.
\end{multline*}
Thus if $T < \varrho$, then
\[
\biggl(\sup_{B_\varrho(y) \cap U} v_p - \inf_{B_\varrho(x) \cap U} v_p\biggr)^{\frac{1}{p - 1}} 
\ge \frac{1 - 2\epsilon}{(1 + \epsilon)^2} \ell^{\frac{1}{p - 1}} \fint_0^\ell (\lambda'(s))^\perp \cdot \nabla u_\infty(\lambda(s)) \, ds.
\]
Since $\epsilon > 0$ was chosen arbitrarily, the second inequality follows.
\end{proof}

Now fix a constant $\gamma \in (0, 1 - \delta)$.
Since we may add arbitrary constants to $v_p$ without changing
the properties used, we may assume that $v_p(-\frac{3}{4}, 0) = \gamma^{p - 1}$ for
every $p \in [2, \infty)$. In view of inequality \eqref{eqn:angle_v} and the
considerations immediately following it,
if $\delta$ and $1 - \sigma$ are sufficiently small and $p$ is
sufficiently large, then
$v_p \le \gamma^{p - 1}$ in $(-1, -\frac{7}{8}] \times [-\frac{1}{2}, \frac{1}{2}]$
and $v_p \ge \gamma^{p - 1}$ in $[-\frac{5}{8}, 1) \times [-\frac{1}{2}, \frac{1}{2}]$.
In particular, the functions
\[
w_p = \frac{\log v_p}{1 - p}
\]
are well-defined and satisfy $w_p \le -\log \gamma$ at least in $[-\frac{5}{8}, 1) \times [-\frac{1}{2}, \frac{1}{2}]$. We observe that
\begin{equation} \label{eqn:imcf-p}
\Delta_{p'} w_p = |\nabla w_p|^{p'}
\end{equation}
wherever $w_p$ is defined. Thus for any $\eta \in C_0^\infty(Q_{1/2})$ with $\eta \ge 0$,
\[
\begin{split}
\int_U \eta^{p'} |\nabla w_p|^{p'} \, dx & = -p' \int_U \eta^{p' - 1} |\nabla w_p|^{p' - 2} \nabla \eta \cdot \nabla w_p \, dx \\
& \le p'\left(\int_U |\nabla \eta|^{p'} \, dx\right)^{\frac{1}{p'}} \left(\int_U \eta^{p'} |\nabla w_p|^{p'} \, dx\right)^{\frac{1}{p}}.
\end{split}
\]
Hence
\[
\int_U \eta^{p'} |\nabla w_p|^{p'} \, dx \le (p')^{p'} \int_U |\nabla \eta|^{p'} \, dx.
\]
This means that $\|\nabla w_p\|_{L^1(K)}$ is uniformly bounded for any precompact set $K \Subset Q_{1/2}$.

We apply Lemma \ref{lem:pointwise} to $\lambda(s) = (s - \frac{15}{16}, 0)$
for $0 \le s \le \frac{27}{16}$.
We already know that
\[
\sup_{B_{1/16}(-15/16, 0)} v_p \le \gamma^{p - 1}.
\]
The first inequality in Lemma \ref{lem:pointwise} then gives the estimate
\[
\inf_{B_{1/16}(3/4, 0)} v_p \le (1 + 2\delta)^{p - 1} + \gamma^{p - 1}
\]
for $p$ sufficiently large.
Using \eqref{eqn:angle_v} again, we conclude that
\[
v_p \le (1 + 2\delta)^{p - 1} + \gamma^{p - 1}
\]
in $Q_{1/2}$. It follows that $w_p$ is uniformly bounded in $L^\infty(Q_{1/2})$.
From \cite[Proposition 2.1]{Moser:15.2}, we then obtain a uniform bound for $\|w_p\|_{W^{1, q}(K)}$ for any
precompact $K \Subset Q_{1/2}$ and any $q < \infty$. Therefore, there exists a
sequence $p_k \to \infty$ such that $w_{p_k} \rightharpoonup w$ weakly in $\bigcap_{q < \infty} W_\loc^{1, q}(Q_{1/2})$
for some function $w \colon Q_{1/2} \to \R$. By the Sobolev embedding theorem
and the Arzel\`a-Ascoli theorem, the convergence is also locally uniform and $w$
is continuous.

\begin{lemma} \label{lem:pointwise2}
Let $\lambda \colon [0, \ell] \to U$ be a solution of the equation
\[
\lambda'(s) = -\frac{\nabla^\perp u_\infty(\lambda(s))}{|\nabla^\perp u_\infty(\lambda(s))|}
\]
for $s \in [0, \ell]$. Let $x = \lambda(0)$ and $y = \lambda(\ell)$.
If $x \in (-1, -\frac{15}{16}] \times [-\frac{3}{8}, \frac{3}{8}]$
and $y \in Q_{1/4}$, then
\[
\inf_{0 \le s \le \ell} \bigl(- \log |\nabla u_\infty(\lambda(s))|\bigr) \le w(y) \le - \log |\nabla u_\infty(y)|.
\]
\end{lemma}

\begin{proof}
Set $A = \sup_{0 \le s \le \ell} |\nabla u_\infty(\lambda(s))|$.
We apply Lemma \ref{lem:pointwise} with $\varrho \le \frac{1}{16}$.
Then
\[
\sup_{B_\varrho(x) \cap U} v_p \le \gamma^{p - 1}.
\]
Let $\epsilon > 0$. It follows that for $p$ sufficiently large,
\[
\inf_{B_\varrho(y)} v_p \le (A+ \epsilon)^{p - 1} + \gamma^{p - 1},
\]
and therefore,
\[
\begin{split}
\sup_{B_\varrho(y)} w_p & \ge \frac{1}{1 - p} \log \left((A+ \epsilon)^{p - 1} + \gamma^{p - 1}\right) \\
& = -\log (A + \epsilon) + \frac{1}{1 - p} \log \left(1 + \left(\frac{\gamma}{A + \epsilon}\right)^{p - 1}\right).
\end{split}
\]
Since $\gamma < 1 - \delta \le A$,
it follows that
\[
\sup_{B_\varrho(y)} w \ge -\log (A + \epsilon).
\]
Since $w$ is continuous, we obtain the first estimate by letting $\varrho \to 0$ and $\epsilon \to 0$.

For the proof of the second inequality, we apply Lemma \ref{lem:pointwise}
to the restriction of $\lambda$ to $[k, \ell]$, where $k \in [0, \ell)$
is chosen such that $\lambda(k) \in Q_{1/2}$ and
\[
\fint_k^\ell |\nabla u_\infty(\lambda(s))| \, ds \ge |\nabla u_\infty(y)| - \epsilon.
\]
Then
\[
\liminf_{p \to \infty} \left(\sup_{B_\varrho(y)} v_p - \inf_{B_\varrho(\lambda(k))} v_p\right)^{\frac{1}{p - 1}} \ge |\nabla u_\infty(y)| - \epsilon.
\]
Note that $\inf_{B_\varrho(\lambda(k))} v_p \ge \gamma^{p - 1} > 0$ by the above observations.
Hence for $p$ large enough,
\[
\sup_{B_\varrho(y)} v_p \ge (|\nabla u_\infty(y)| - 2\epsilon)^{p - 1}
\]
and
\[
\inf_{B_\varrho(y)} w_p \le - \log (|\nabla u_\infty(y)| - 2\epsilon).
\]
The same inequality follows for $w$.
Again we conclude the proof by letting $\varrho \to 0$ and $\epsilon \to 0$.
\end{proof}

Under the assumptions of Proposition \ref{prp:local}, the level sets
$L_t$ of $u_\infty$ can be parametrised by curves $\lambda$ as in
Lemma \ref{lem:pointwise2}. Considering the monotonicity of
$|\nabla u_\infty|$ along these level sets, it follows that
\[
\sup_{0 \le s \le \ell} |\nabla u_\infty(\lambda(s))| = |\nabla u_\infty(\lambda(\ell))|
\]
if $\lambda(\ell) \in M$. Furthermore, if $\delta$ is sufficiently small,
then any $L_t$ that intersects $Q_{1/4}$ will also intersect
$(-1, -\frac{15}{16}] \times [-\frac{3}{8}, \frac{3}{8}]$. Thus Lemma \ref{lem:pointwise2} implies
that $w = - \log |\nabla u_\infty|$ in $Q_{1/4} \cap M$, which is statement \ref{itm:2}.

Next we have a closer look at equation \eqref{eqn:imcf-p} and study what it means for the limit.
We define $F_p = |\nabla w_p|^{p' - 2} \nabla w_p$, so we can write
\begin{equation} \label{eqn:F_p}
\div F_p = |\nabla w_p|^{p'}.
\end{equation}
For any $q < \infty$ and any precompact set $K \Subset Q_{1/2}$, we have the inequality
\[
\left(\int_K |F_p|^q \, dx\right)^{\frac{1}{q}} \le |K|^{\frac{p - q}{pq}} \left(\int_K |F_p|^p \, dx\right)^{\frac{1}{p}} = |K|^{\frac{p - q}{pq}} \left(\int_K |\nabla w_p|^{p'} \, dx\right)^{\frac{1}{p}}.
\]
Hence
\[
\limsup_{p \to \infty} \left(\int_K |F_p|^q \, dx\right)^{\frac{1}{q}} \le |K|^{\frac{1}{q}}.
\]
Therefore, we may assume that $F_{p_k} \rightharpoonup F$ weakly in $L_\loc^q(Q_{1/2}; \R^2)$
for every $q < \infty$. Moreover,
\[
\|F\|_{L^\infty(K)} = \lim_{q \to \infty} \|F\|_{L^q(K)} \le 1.
\]
That is, we have the inequality $|F| \le 1$ almost everywhere in $Q_{1/2}$.

We also observe that by the definition of $w_p$,
\[
F_p = (p - 1)^{-\frac{1}{p - 1}} e^{w_p} \nabla^\perp u_p.
\]
As $\nabla u_{p_k} \rightharpoonup \nabla u_\infty$ weakly in $L_\loc^q(U)$ and $w_{p_k} \to w$ locally
uniformly in $Q_{1/2}$, this implies that
\[
F = e^w \nabla^\perp u_\infty.
\]
It follows that $w \le -\log |\nabla u_\infty|$.

From arguments developed for the inverse mean curvature flow \cite[p.~82]{Moser:07.1}, it follows that
\[
\int_U \eta |\nabla w_p|^{p'} \, dx \to \int_U \eta |\nabla w| \, dx
\]
for any $\eta \in C_0^\infty(Q_{1/2})$. Hence \eqref{eqn:F_p} gives rise to
\[
\div F = |\nabla w|
\]
in $Q_{1/2}$, which amounts to statement \ref{itm:3}.
Testing this equation with $\eta e^{-w}$, we find that
\[
\begin{split}
\int_U \eta e^{-w} |\nabla w| \, dx & = -\int_U e^{-w} (\nabla \eta - \eta \nabla w) \cdot F \, dx \\
& = -\int_U (\nabla \eta - \eta \nabla w) \cdot \nabla^\perp u_\infty \, dx \\
& = \int_U \eta \nabla w \cdot \nabla u_\infty^\perp \, dx
\end{split}
\]
for any $\eta \in C_0^\infty(Q_{1/2})$. Therefore,
\[
e^{-w} |\nabla w| = \nabla w \cdot \nabla^\perp u_\infty
\]
almost everywhere in $Q_{1/2}$. Thus we have proved statement \ref{itm:1} as well.

\section{Weak solutions of the inverse mean curvature flow} \label{sct:imcf}

The study of weak solutions of the inverse mean curvature flow
goes back to a seminal paper of Huisken and Ilmanen \cite{Huisken-Ilmanen:01},
where they are defined in terms of a variational condition.
Equation \eqref{eqn:imcf} is not actually variational, but Huisken and Ilmanen
get around that problem by asking that a function $w$ minimise a functional
depending on $w$ itself. They work with local Lipschitz functions,
but the same ideas make sense for functions
$w \in W_\loc^{1, 1}(\Omega) \cap L_\loc^\infty(\Omega)$.
According to their definition, $w$ is a weak solution of \eqref{eqn:imcf} if 
\begin{equation} \label{eqn:Huisken-Ilmanen}
\int_K (1 + w)|\nabla w| \, dx \le \int_K \left(|\nabla \tilde{w}| + \tilde{w} |\nabla w|\right) \, dx
\end{equation}
for every precompact set $K \Subset \Omega$ and every
$\tilde{w} \in W_\loc^{1, 1}(\Omega) \cap L_\loc^\infty(\Omega)$ with
$\tilde{w} = w$ in $\Omega \setminus K$.

The same paper also proves a number of properties implied
by this condition, including a comparison principle, a resulting uniqueness
theorem, and a minimising hull property for the sublevel sets. Despite
some technical differences, many of the underlying ideas will apply to
the situation discussed in the introduction, too. It is therefore
useful to know that the condition from Definition \ref{def:imcf}
implies inequality \eqref{eqn:Huisken-Ilmanen}, provided that $w \in L_\loc^\infty(\Omega)$.

\begin{proposition}
Suppose that $w \in W_\loc^{1, 1}(\Omega) \cap L_\loc^\infty(\Omega)$ is a weak
solution of \eqref{eqn:imcf} in the sense of Definition \ref{def:imcf}.
Then for every $K \Subset \Omega$ and every
$\tilde{w} \in W_\loc^{1, 1}(\Omega) \cap L^\infty(\Omega)$ with
$\tilde{w} = w$ in $\Omega \setminus K$, inequality \eqref{eqn:Huisken-Ilmanen}
holds true.
\end{proposition}

\begin{proof}
Let $F \in L^\infty(\Omega; \R^n)$ be a vector field satisfying
$|F| \le 1$ and $F \cdot \nabla w = |\nabla w|$ almost everywhere in
$\Omega$ and
\[
\div F = |\nabla w|
\]
weakly. Then for almost every $x \in \Omega$, either $\nabla w(x) = 0$ or
\[
F(x) = \frac{\nabla w(x)}{|\nabla w(x)|}.
\]

Let $K \Subset \Omega$ be a precompact set and let
$\tilde{w} \in W_\loc^{1, 1}(\Omega) \cap L_\loc^\infty(\Omega)$ with
$\tilde{w} = w$ in $\Omega \setminus K$. Then
\[
|\nabla \tilde{w}| \ge F \cdot \nabla \tilde{w} = |\nabla w| + F \cdot (\nabla \tilde{w} - \nabla w)
\]
almost everywhere. It follows that
\[
\int_K |\nabla \tilde{w}| \, dx \ge \int_K |\nabla w| \, dx - \int_K (\tilde{w} - w) |\nabla w| \, dx.
\]
Rearranging the terms, we obtain inequality \eqref{eqn:Huisken-Ilmanen}.
\end{proof}

\def\cprime{$'$}
\providecommand{\bysame}{\leavevmode\hbox to3em{\hrulefill}\thinspace}
\providecommand{\MR}{\relax\ifhmode\unskip\space\fi MR }
\providecommand{\MRhref}[2]{%
  \href{http://www.ams.org/mathscinet-getitem?mr=#1}{#2}
}
\providecommand{\href}[2]{#2}

\end{document}